\documentclass[reqno, oneside]{amsart}
\usepackage{style}

\begin{document}
\title[Truncated Toeplitz operators]{Fredholmness and compactness of truncated Toeplitz and Hankel operators}
\author{R.~V.~Bessonov}

\address{St.Petersburg State University ({\normalfont \hbox{7-9}, Universitetskaya nab., 199034, St.Petersburg,  Russia}), St.Petersburg Department of Steklov Mathematical Institute of Russian Academy of Science ({\normalfont 27, Fon\-tan\-ka, 191023, St.Petersburg, Russia}), and School of Mathematical Sciences, Tel Aviv University ({\normalfont 69978, Tel Aviv, Israel})}
\email{bessonov@pdmi.ras.ru}
\thanks{This work is partially supported by RFBR grants 12-01-31492, 14-01-00748, by ISF grant~94/11, by JSC ``Gazprom Neft'' and by the Chebyshev Laboratory (Department of Mathematics and Mechanics, St. Petersburg State University) under RF Government grant 11.G34.31.0026}

\begin{abstract}
We prove the spectral mapping theorem $\sigma_e(A_\phi) = \phi(\sigma_e(A_z))$ for the Fredholm spectrum of a truncated Toeplitz operator $A_\phi$ with symbol $\phi$ in the Sarason algebra $\sal$ acting on a coinvariant subspace $\Kth$ of the Hardy space $H^2$. Our second result says that a truncated Hankel operator on the subspace $\Kth$ generated by a one-component inner function $\theta$ is compact if and only if it has a continuous symbol. We also suppose a description of truncated Toeplitz and Hankel operators in Schatten classes $S^p$.
\end{abstract}

\maketitle

\section{Introduction}\label{s1}
Truncated Toeplitz and Hankel operators are compressions of the standard Toep\-litz and Hankel operators on the Hardy space $H^2$ to its coinvariant subspaces $\Kth$. More precisely, consider the shift operator $S: f \mapsto zf$ on the Hardy space~$H^2$ in the open unit disk $\D$ of the complex plane $\C$. By A. Beurling theorem, invariant subspaces of $S$ have the form $\theta H^2$, where $\theta$ is an inner function in $\D$. Accordingly, the subspaces $\Kth = H^2\ominus\theta H^2$ are invariant under the backward shift operator $S^*: f \mapsto \frac{f - f(0)}{z}$ on $H^2$. Fix an inner function $\theta$ and consider $H^2$ and $\Kth$ as subspaces of the space $L^2 = L^2(\T)$ on the unit circle~$\T$. Let $P_\theta$ and $P_{\bar{\theta}}$ denote the orthogonal projections in $L^2$ to the subspaces $\Kth$ and $\ov{z\Kth} = \{f \in L^2: \ov{zf} \in \Kth\}$, correspondingly. Take $\phi \in L^2$ and define the operators 
$$
A_\phi : f\mapsto P_{\theta} (\phi f), \qquad \Gamma_\phi : f\mapsto P_{\bar\theta} (\phi f), 
$$
on the dense subset $\Kth \cap L^\infty$ of the space $\Kth$. The operator $A_\phi: \Kth \to \Kth$ is called the truncated Toeplitz operator; $\Gamma_\phi: \Kth \to \ov{z\Kth}$ is the truncated Hankel operator with symbol $\phi$.
For recent results on truncated Toeplitz operators see survey \cite{GR13}. The present paper is closely related to  \cite{BBK11}, \cite{BCFMT}, \cite{BR15}, \cite{Roch87}, \cite{Sar07}. Let us now discuss in details two results formulated in the abstract.

\subsection*{A spectral mapping theorem} The restricted shift operator $\Sth = A_z$ on $\Kth$ is the simplest non-trivial case of Sz.-Nagy--Foias model \cite{NF10} for Hilbert space contractions. Textbook \cite{Nik86} by N. K. Nikolski is devoted to the study of different properties of this operator. 
Defining $\phi(\Sth) = A_\phi$ for bounded analytic functions $\phi \in H^\infty$, we obtain the $H^\infty$-functional calculus for the operator $\Sth$. Let $\sigma(\theta)$ denote the spectrum of the inner function 
$\theta$, that is, the set of all points $\xi$ in the closed unit disk such that $\liminf_{|z|< 1, \, z \to \xi} |\theta(z)| = 0$. Then (see Section III.3 in \cite{Nik86}) for every function $\phi \in H^\infty$ we have
\begin{equation}\label{eq2}
\sigma(A_\phi) = \phi(\sigma(\theta)),
\end{equation}
where $\sigma(A_\phi)$ denotes the spectrum of the operator $A_\phi$. Formula \eqref{eq2} can be regarded as a spectral mapping theorem for the $H^\infty$-functional calculus of the operator $\Sth$. Indeed, for every $\phi \in H^\infty$ we have $\sigma(\phi(\Sth)) = \phi(\sigma(\theta)) = \phi(\sigma(S_\theta))$. 

\medskip

Proof of formula $(1)$ relies on the celebrated Corona theorem for the algebra~$H^\infty$ of bounded analytic functions in the unit disk $\D$. On the other hand, by much more elementary technique one can show that for every continuous function $\phi$ on the unit circle $\T$ we have
\begin{equation}\label{eq5}
\sigma_{e}(A_\phi) = \phi(\sigma(\theta) \cap \T),
\end{equation} 
where $\sigma_{e}(A_\phi)$ is the Fredholm (or essential) spectrum of $A_\phi$, see Section V.4 in \cite{Nik86} or \cite{GRW12}. Denote by $\ct$ the set of all continuous functions on $\T$. 
D. Sarason \cite{Sar67} proved that the set 
$$
\sal = \{\phi_1 + \phi_2, \;\; \phi_1 \in \ct, \; \phi_2 \in H^{\infty}\}
$$
is the closed subalgebra of $L^{\infty}$. We compute the essential spectrum of truncated Toeplitz operators with symbols in $\sal$.

\begin{Thm}\label{t1}
Let $\theta$ be an inner function, and let $A_\phi$ be the truncated Toeplitz operator on $\Kth$ with symbol $\phi \in \sal$. Then $\sigma_{e}(A_\phi) = \phi(\sigma_{e}(A_z)) = \phi(\sigma(\theta) \cap \T)$. 
\end{Thm}
The main step in the proof of Theorem \ref{t1} is an application of the corona theorem for the Sarason algebra $\sal$ obtained  in 2007 by R. Mortini and B. Wick \cite{MW10}. 

\medskip

Let $B(\Kth)/S^{\infty}(\Kth)$ be the Calkin algebra of all bounded operators on $\Kth$ modulo compact operators. We will denote its elements by $[T]$. Note that for $T \in B(\Kth)$ we have $\sigma_e(T) = \sigma([T])$ by the definition of the essential spectrum. It not difficult to check that the mapping $\phi \mapsto [A_\phi]$ is the contractive homomorphism from $\sal$ to $B(\Kth)/S^{\infty}(\Kth)$. Hence if we define the $\sal$-functional calculus of $[\Sth]$ by $\phi([\Sth]) = [A_\phi]$ for $\phi \in \sal$, then Theorem \ref{t1} becomes the spectral mapping theorem for this calculus. 

\medskip

\subsection*{Compact truncated Hankel operators} 
The classical result by Hartman \cite{Hart58} says that a Hankel operator $H_\phi:H^2\to\ov{zH^2}$ is compact if and only if it has a continuous symbol. Its modern proof (see, e.g., Section 1.5 in \cite{PeBook}) is based on the Nehari theorem $\|H_\phi\| = \dist_{L^\infty}(\phi, H^\infty)$ and the fact that
for every function $\phi \in \ct$ we have
\begin{equation}\label{eq3}
\dist_{L^\infty}(\phi, H^\infty) = \dist_{L^\infty}(\phi, H^\infty \cap \ct).  
\end{equation}  
Formula \eqref{eq3} follows easily from properties of the Poisson kernel. Another proof is due to D. Sarason. It uses a corollary of brothers Riesz theorem on analytic measures, namely, the duality relation $\bigl(\ct/\mathcal{A}\bigr)^{**} = L^{\infty}/H^{\infty}$, where $\mathcal{A} = H^\infty \cap \ct$ is the disk algebra. Both proofs can be found in Section VII of \cite{Koos98}. As we will see, the duality approach works for the following analogue of \eqref{eq3}: for every $\phi \in \ct$ we have
\begin{equation}\label{eq6}
\dist_{L^\infty}(\phi, \ft) = \dist_{L^\infty}(\phi, \ft \cap \ct), 
\end{equation}
where $\ft$ denotes the closure of the set $\ov{\theta H^\infty} + H^{\infty} =\{\ov{\theta f_1} + f_2, \; f_{1,2} \in H^\infty\}$ in $w^*$-topology of the space $L^{\infty}$. It worth be mentioned that the elementary argument based on a convolution with the Poisson kernel gives \eqref{eq6} only in the case where $\theta$ is a finite Blaschke product, see Section \ref{s3}. 

\medskip

An inner function $\theta$ is said to be one-component if the set $\{z \in \D: \; |\theta(z)| < \eps \}$ is connected for a positive number~$\eps < 1$. It follows from the results of \cite{BBK11} that truncated Hankel operators acting on the space $\Kth$ generated by a one-component inner function $\theta$ satisfy the two-sided estimate $\|\Gamma_\phi\| \asymp \dist_{L^\infty}(\phi, \fts)$ of Nehari type. Using this estimate and formula~\eqref{eq6} for the inner function $\theta^2$, we obtain the following theorem. 
\begin{Thm}\label{t2}
Let $\theta$ be a one-component inner function. Then a truncated Hankel operator $\Gamma: \Kth \to \ov{\Kth}$ is compact if and only if $\Gamma = \Gamma_{\phi}$ for some $\phi \in \ct$. Moreover, one can choose $\phi \in \ct$ so that 
$\|\Gamma\| \le \|\phi\|_{L^{\infty}} \le c_\theta\|\Gamma\|$, where the constant $c_\theta$ depends only 
on~$\theta$. 
\end{Thm} 
In Section \ref{s4} we combine Theorem \ref{t2} with some other characterizations of compact truncated Hankel operators. We also present a conjecture on truncated Hankel operators in Schatten classes $S^p(\Kth)$. We expect that their description is possible in terms of a mean oscillation with respect to the Clark measure of the inner function~$\theta^2$.     

\section{Proof of Theorem \ref{t1}}\label{s2}
Let us recall some definitions. A bounded operator $T$ on the Hilbert space $H$ is called Fredholm if its range $\ran T$ is the closed subspace of $H$, $\dim \ker T < \infty$, and $\dim \ker T^* < \infty$. A basic Fredholm theory says that $T$ is Fredholm if and only if there are bounded operators $L$, $R$ and compact operators $K_L$, $K_R$ on $H$ such that $LT=I+K_L$, $TR=I+K_R$, where~$I$ denotes the identical operator on $H$. The essential (Fredholm) spectrum of $T$, $\sigma_e(T)$, is the set of all $\lambda \in \C$ such that the operator $T - \lambda I$ is not Fredholm. Define the continuous spectrum $\sigma_c(T)$ of $T$ as the set of all $\lambda \in \C$ such that there exists a non-compact sequence $\{x_n\} \subset H$ for which $\lim_n\|(T - \lambda I)x_n\| = 0$. It is clear that $\sigma_c(T) \subset \sigma_e(T)$. Moreover, one can check that $\sigma_e(T) = \sigma_c(T) \cup \ov{\sigma_c(T^*)}$.

\medskip

We will prove that $\sigma_e(A_\phi) \subset \phi (\sigma(\theta) \cap \T) \subset \sigma_c(A_\phi)$ for every $\phi \in \sal$. Since functions $\phi \in \sal$ are not defined everywhere on $\T$, we need a definition of the image $\phi (\sigma(\theta) \cap \T)$. Let $m$ be the normalized Lebesgue measure on the unit circle~$\T$. Put
\begin{equation}\notag
\phi (\sigma(\theta) \cap \T) = \bigr\{\zeta \in \C: \; \liminf_{z \in \D, \, |z| \to 1}(|\hat{\phi}(z) - \zeta| + |\theta(z)|) = 0\bigl\},
\end{equation}  
where $\hat \phi$ denotes the Poisson transform of $\phi$, 
$$\hat{\phi}(z) = \int_{\T}\phi(\xi) \frac{1 - |z|^2}{|1 - \bar \xi z|^2}\,dm(\xi), \quad z \in \D.$$ 
For a function $\phi \in \ct$ thus defined set $\phi (\sigma(\theta) \cap \T)$ coincides with the usual image of $\sigma(\theta) \cap \T$ by $\phi$. 

\medskip

Proof of Theorem \ref{t1} is based on the following result \cite{MW10}. 
\begin{NMT}(R.Mortini, B. Wick)
Let $f_1, \ldots, f_n$ be a family of functions in $\sal$ such that for some $\eps>0$ and $r < 1$ we have 
$$|\hat{f}_1(z)| + \ldots + |\hat{f}_n(z)| > \eps, \quad r \le |z| < 1.$$
Then there exist functions $g_1, \ldots, g_n$ in $\sal$ such that 
$f_1 g_1 + \ldots + f_n g_n = 1$ almost everywhere on the unit circle $\T$. 
\end{NMT} 
The corresponding statement in \cite{MW10} (Theorem 1.1) is given in slightly different terms. However, that statement can be easily reduced to the above formulation by noting that $f \in \sal$ is invertible if and only if there exists $r < 1$ such that $|\hat f(z)| > \eps$ for some $\eps>0$ and all $z \in \D$ with $r \le |z| < 1$. Alternatively, one can check that the proof in \cite{MW10} works for our reformulation without any essential changes.   
 
\begin{Lem}\label{l1}
Let $\theta$ be an inner function and let $\phi \in \sal$. Then for the truncated Toeplitz operator $A_\phi: \Kth \to \Kth$ we have $\phi(\sigma(\theta) \cap \T) \subset \sigma_{c}(A_\phi)$.
\end{Lem} 
\beginpf The proof of lemma is a modification of arguments in Section III.3 of \cite{Nik86}. Let $\phi \in \sal$ and let $\zeta \in \phi(\sigma(\theta) \cap \T)$. Then there exists a sequence $\{\lambda_n\} \subset \D$ such that 
$$\lim_{n \to \infty}(|\hat{\phi}(\lambda_n) - \zeta| + |\theta(\lambda_n)|) = 0.$$
We can assume that $\lambda_n$ converge to a point $\lambda_{\infty} \in \T$ and that $|\theta(\lambda_n)| < \frac{1}{2}$ for all $n$. For $\lambda \in \D$ denote by $\tilde k_{\lambda}$ the function 
$\tilde k_{\lambda} = \frac{\theta - \theta(\lambda)}{z - \lambda}$. It is easy to check that $\tilde k_{\lambda} \in \Kth$ and $\|\tilde k_{\lambda}\|^2 = \frac{1-|\theta(\lambda)|^2}{1 - |\lambda|^2}$. We claim that 
\begin{equation}\label{eq14}
\lim_{n \to \infty} \frac{\|(A_\phi - \zeta I)\tilde k_{\lambda_n}\|}{\|\tilde k_{\lambda_n}\|} = 0.
\end{equation}
Consider functions $\phi_1 \in \ct$, $\phi_2 \in H^{\infty}$ such that $\phi = \phi_1 + \phi_2$. Put $\zeta_1 = \phi_1(\lambda_{\infty})$ and $\zeta_2 = \zeta - \zeta_1$. We have
\begin{equation}\label{eq15}
\frac{\|(A_{\phi_1} - \zeta_1 I)\tilde k_{\lambda_n}\|^2}{\|\tilde k_{\lambda_n}\|^2} \le 
\frac{\|(\phi_1 - \zeta_1)\tilde k_{\lambda_n}\|^2}{\|\tilde k_{\lambda_n}\|^2} \le 
8\int_{\T} |\phi_1(z) - \zeta_1|^2 \frac{1 - |\lambda_n|^2}{|z - \lambda_n|^2} \, dm(z). 
\end{equation}
Since the function $|\phi_1 - \zeta_1|^2$ is continuous on $\T$ and vanishes at $\lambda_\infty = \lim \lambda_n$, the last integral tends to zero as $n \to \infty$. Next, for every $\lambda_n$ the function $\frac{\phi_2 - \phi_2(\lambda_n)}{z - \lambda_n}$ belongs to $H^2$, hence
$$
A_{\phi_2} \tilde k_{\lambda_n} = P_{\theta} \left( \phi_2 \frac{\theta - \theta(\lambda_n)}{z - \lambda_n}\right) = -\theta(\lambda_n) P_{\theta} \left( \frac{\phi_2 - \phi_2(\lambda_n)}{z - \lambda_n}\right) + \phi_2(\lambda_n)\frac{\theta - \theta(\lambda_n)}{z - \lambda_n}.
$$ 
Using the estimate $\bigl\|\frac{\phi_2 - \phi_2(\lambda_n)}{z - \lambda_n}\bigr\|^2 \le \frac{4\|\phi_2\|_{L^\infty}^{2}}{1-|\lambda_n|^2}$, we obtain
\begin{equation}\label{eq16}
\frac{\|(A_{\phi_2} - \zeta_2 I)\tilde k_{\lambda_n}\|^2}{\|\tilde k_{\lambda_n}\|^2} \le 8|\theta(\lambda_n)|^2 \|\phi_2\|_{L^\infty}^{2}+|\phi_2(\lambda_n) - \zeta_2|^2. 
\end{equation}
Both summands in the right hand side tend to zero as $n \to \infty$. Combining estimates~\eqref{eq15} and \eqref{eq16}, we see that \eqref{eq14} holds. Since $\lim|\theta(\lambda_n)| =0$, the sequence $\tilde k_{\lambda_n}/\|\tilde k_{\lambda_n}\|$ is non-compact (use the fact that the sequence $\sqrt{1 - |\lambda_n|^2}/(1-\bar\lambda_n z)$ is non-compact in $H^2$). It follows that $\phi(\sigma(\theta) \cap \T) \subset \sigma_{c}(A_\phi)$. \qed

\begin{Lem}\label{l3}
Let $\theta$ be an inner function. The mapping $\phi \mapsto [A_\phi]$ is the contractive homomorphism from $\sal$ to $B(\Kth)/S^{\infty}(\Kth)$.
\end{Lem}
\beginpf Clearly, the mapping $\phi \mapsto [A_\phi]$ is linear and $\|[A_\phi]\| \le \|A_\phi\| \le \|\phi\|_{L^\infty}$. We have $A_{\phi_1} A_{\phi_2} = A_{\phi_1 \phi_2}$ for all $\phi_1$, $\phi_2$ in $H^{\infty}$, see Section III.2 in \cite{Nik86}. Also, it was observed in \cite{GRW12} that if $\phi_1 \in \ct$ and $\phi_2 \in L^{\infty}$, then $A_{\phi_1} A_{\phi_2} = A_{\phi_1 \phi_2} + K_1$ and $A_{\phi_2} A_{\phi_1} = A_{\phi_1 \phi_2} + K_2$, where $K_1$, $K_2$ are compact operators on $\Kth$. This yields the statement of the lemma. \qed
\begin{Lem}\label{l2}
Let $\theta$ be an inner function and let $\phi \in \sal$. Then for the truncated Toeplitz operator $A_\phi: \Kth \to \Kth$ we have $\sigma_{e}(A_\phi) \subset \phi(\sigma(\theta) \cap \T)$.  
\end{Lem}
\beginpf Take a point $\zeta \in \C\setminus\phi(\sigma(\theta) \cap \T)$. Let us show that the operator $A_\phi - \zeta I$ on $\Kth$ is Fredholm. We have
\begin{equation}\label{eq4}
|\hat \phi(z) - \zeta| + |\theta(z)| > \eps, \qquad z \in U_\delta \cap \D,
\end{equation}
for some neighbourhood $U_\delta$ of the closed set $\sigma(\theta) \cap \T$. From the definition of $\sigma(\theta)$ we see that estimate \eqref{eq4} holds for all $z$ with $|z| \ge r$, where $r<1$ is a positive number depending on $\theta$ and $\eps$. Thus, we can apply the corona theorem for the algebra $\sal$ and find functions $g_1$, $g_2$ in $\sal$ such that
\begin{equation}\label{eq23}
(\phi - \zeta)g_1 + \theta g_2 = 1
\end{equation}
almost everywhere on the unit circle $\T$. Note that $A_\theta$ is the zero operator on $\Kth$. Using 
equation \eqref{eq23} and Lemma \ref{l3}, we see that  
$$A_{g_1} (A_\phi - \zeta I) = I + K_L, \qquad (A_\phi - \zeta I) A_{g_1} = I + K_R,$$
for some compact operators $K_L$, $K_R$  on $\Kth$. It follows that $\zeta \in \C \setminus \sigma_e(A_\phi)$, as required. \qed

\medskip

\noindent{\bf Proof of Theorem \ref{t1}.} By Lemma \ref{l1} and Lemma \ref{l2} we have the following chain of inclusions: $\sigma_e(A_\phi) \subset \phi (\sigma(\theta) \cap \T) \subset \sigma_c(A_\phi)$. Since $\sigma_c(T) \subset \sigma_e(T)$ for every bounded operator $T$, we have 
$\sigma_e(A_\phi) = \phi (\sigma(\theta) \cap \T)$. In particular, for $\phi = z$ we get $\sigma_e(A_z) = \sigma(\theta) \cap \T$ and hence $\sigma_e(A_\phi) = \phi (\sigma_e(A_z))$. 
The theorem is proved. \qed

\medskip

Next proposition describes the kernel of the homomorphism $\phi \mapsto A_\phi$ from $\sal$ to $B(\Kth)/S^{\infty}(\Kth)$. It is an extension of the Sarason theorem on compact truncated Toeplitz operators with symbols in $H^\infty$, see Section VIII.3 in \cite{Nik86}.  

\begin{Prop}\label{p2}
Let $\theta$ be an inner function and let $\phi \in \sal$. Then the truncated Toeplitz operator $A_\phi: \Kth \to \Kth$ is compact if and only if $\phi \in \theta \ct + \theta H^{\infty}$.
\end{Prop}
\beginpf Let $P_{-}$ denote the orthogonal projector on $L^2$ to the subspace $\ov{z H^2}$. Consider the operator 
$T: g \mapsto \theta P_{-} (\bar \theta g) - P_{-} g$ on $L^2$. It is easy to see that $Tg = 0$ for $g \in \ov{zH^2} \oplus \theta H^2$ and $Tg = g$ for $g\in \Kth$. Hence, we have $T = P_\theta$. Define the Hankel operator 
$H_{\psi}: H^2 \to \ov{zH^2}$ with symbol $\psi \in L^\infty$ by
$H_{\psi}: f \mapsto P_{-}(\psi f)$. By Hartman's theorem, $H_{\psi}$ is compact if and only if $\psi \in \sal$, see Section 1.5 in \cite{PeBook}. For every function $f \in \Kth$ we have 
\begin{equation}\label{eq12}
A_\phi f = P_\theta (\phi f) = \theta P_{-} (\bar \theta \phi f) - P_{-} (\phi f) = \theta H_{\bar \theta \phi} f - H_{\phi} f.  
\end{equation}
Assume that $\phi \in \sal$ is such that $A_\phi \in S^{\infty}(\Kth)$. Consider a weakly convergent sequence 
$g_n \in H^2$ with zero limit. Put $f_n = P_{\theta} g_n$. From formula \eqref{eq12} we see that $\lim\|H_{\bar \theta \phi} f_n\| = 0$ because the operators $A_\phi$, $H_\phi$ are compact. Consider the sequence $\theta h_n = g_n - f_n$ of functions in $\theta H^2$. Note that $h_n \in H^2$ converge weakly to zero. Since $\phi \in \sal$, we have 
$\lim\|H_{\bar \theta \phi} \theta h_n\| = \lim\|H_{\phi} h_n\| = 0$. 
It follows that $\lim\|H_{\bar \theta \phi} g_n\| = 0$ for every weakly convergent to zero sequence $\{g_n\} \subset H^2$. Hence the operator $H_{\bar \theta \phi}: H^2 \to \ov{z H^2}$ is compact and therefore $\phi \in \theta \ct + \theta H^\infty$. 

Now let $\phi = \theta \phi_1 + \theta \phi_2$, where $\phi_1 \in \ct$ and $\phi_2 \in H^\infty$. Then we have $A_{\phi} = A_{\theta \phi_1}$. Moreover, from Lemma \ref{l3} we see that $A_{\theta \phi_1} = A_{\theta} A_{\phi_1} + K$ for a compact operator $K$ on $\Kth$. Since $A_{\theta} = 0$, the operator $A_{\theta \phi_1}$ is compact. \qed 

\medskip

\noindent{\bf Remark.} Proposition \ref{p2} can be reformulated in the following way: the truncated Hankel operator $\Gamma_\phi: \Kth \to \ov{z\Kth}$ with symbol $\phi \in \bar{\theta}(\sal)$ is compact if and only if $\phi \in \sal$ (see Lemma \ref{l14}). In general, the assumption 
$\phi \in \bar{\theta}(\sal)$ can not be omitted: one can construct an inner function $\theta$ and a rank-one truncated Hankel operator on $\Kth$ which has no bounded symbols \cite{BCFMT}. In Theorem~\ref{t2} we prove that if $\theta$ is a one-component inner function, then 
{\it every} compact truncated Hankel operator on $\Kth$ has a continuous symbol.

\section{Proof of Theorem \ref{t2}}\label{s3}
The proof of Theorem \ref{t2} splits into series of lemmas. The main analytic ingredient is due to T. Wolff \cite{Wolff82}.
\begin{NMT}[T. Wolff]
Denote $\qc = \ov{(\sal)} \cap (\sal)$ and $\qa = \qc \cap H^\infty$. For every $f \in L^{\infty}$ there exists an outer function $g \in \qa$ such that $gf \in \qc$.
\end{NMT}
Invariant subspaces of the backward shift operator $S^*$ on the Hardy space $H^1$ have the form $\Kl = H^1 \cap \bar z \theta\ov{H^1}$, where $\theta$ is an inner function. It follows that from the definition that a function $f \in L^1$ belongs to $\Kl$ if $\int_{\T} f z h \,dm = 0$ and $\int_{\T} f \ov{\theta h} \,dm = 0$ for all $h \in H^\infty$, where $m$ denotes the normalized Lebesgue measure on the unit circle $\T$. Recall that $\ft$ is the closure of the set $\ov{\theta H^\infty} + H^{\infty} =\{\ov{\theta f_1} + f_2, \; f_{1,2} \in H^\infty\}$ in 
the $w^*$-topology of the space $L^{\infty}$. The following lemma is close to Remark~4.1 in  \cite{Singh94}.
\begin{Lem}\label{l4}
Let $\theta$ be an inner function. Then $\bigl(\ct/(\ft\cap\ct)\bigr)^{*} = \Kl \cap z H^1$. In particular, $\Kl \cap zH^1$ is closed in the $w^*$-topology of the space $zH^1$.
\end{Lem}
\beginpf Let $\Phi$ be a continuous linear functional on $\ct/(\ft\cap\ct)$. By the Riesz-Markov theorem, there exists a complex measure $\mu$ on the unit circle $\T$ such that 
$$
\Phi(\phi) = \int_{\T} \phi \,d\mu, \quad \phi \in \ct,
$$ 
and $\Phi(\phi) = 0$ for all $\phi \in \ft \cap \ct$. This measure $\mu$ can be chosen so that its total variation equals $\|\Phi\|$. For all integers $k \ge 0$ we have $\Phi(z^k) = 0$. Hence $\mu = h \,dm$ for a function $h \in z H^1$ by the brothers Riesz theorem. Using Wolff's theorem, find an outer function $g$ such that $\theta g \in \qc$. Note that $\theta g \in H^\infty$ and therefore $\theta g \in \qa$. For every analytic polynomial $p$ we have $\theta g p \in \qa$ because $\qa$ is an algebra. Since $\qa \subset\ov{\sal}$, one can find functions $\phi \in \ct$, $\psi \in H^{\infty}$ such that 
$\ov{\phi} + \ov{\psi} = \theta g p$. By the construction, we have $\phi \in \ft \cap \ct$. Hence,
$$
\int_{\T} \ov{\theta g p} \, h\,dm = \int_{\T} (\phi + \psi) h\,dm = \int_{\T} \phi h\,dm =\Phi(\phi) = 0.  
$$
Since $g$ is outer, functions of the form $gp$, where $p$ is an analytic polynomial, form the dense subset of $H^{\infty}$ in the $w^*$-topology of $L^\infty$. We see that 
\begin{equation}\notag
\int_{\T} \ov{f} \cdot \ov{\theta} h\,dm = 0  
\end{equation}
for all $f \in H^{\infty}$. It follows that $\ov{\theta} h \in \ov{z H^1}$ and hence $h$ lies in $zH^1 \cap  \bar{z} \theta \ov{H^1} = \Kl \cap zH^1$. By our choice of the measure $\mu = h\,dm$, we have $\|h\|_{L^1} = \|\Phi\|$. Thus, we have the isometric inclusion $(\ct/\ft)^{*} \subset \Kl \cap z H^1$. To obtain the inverse inclusion, observe that for $h \in \Kl \cap z H^1$ the linear functional $\Phi_h: \phi \mapsto \int_{\T} \phi h \, dm$ is continuous on $\ct$ and vanishes on $\ft \cap \ct$. The latter also shows that $\Kl \cap zH^1$ is the annihilator of the subset $\ft \cap \ct  + \mathcal{A}$ of $\ct/\mathcal{A}$. Hence $\Kl \cap zH^1$ is closed in $w^*$-topology of the space~$zH^1$ generated by $\ct/\mathcal{A}$. \qed

\begin{Lem}\label{l5}
Let $\theta$ be an inner function. Then for every $\phi \in \ct$ we have 
$$\dist_{L^\infty}(\phi, \ft) = \dist_{L^\infty}(\phi, \ft \cap \ct).$$
\end{Lem}
\beginpf For every inner function $\theta$ we have $(\Kl \cap zH^1)^* = L^{\infty}/\ft$ by a standard duality argument. Indeed, this follows from the simple fact that $\ft$ is the annihilator of $\Kl \cap zH^1$ in $L^{\infty}$. Hence $\bigl(\ct/(\ft \cap \ct)\bigr)^{**} = L^{\infty}/\ft$ by Lemma \ref{l4}. Consider the canonical embedding 
$$\ct/(\ft \cap \ct) \into L^{\infty}/\ft$$
which sends the factor-class $\phi + \ft \cap \ct$ from $\ct/(\ft \cap \ct)$ to the factor-class $\phi + \ft$ in $L^{\infty}/\ft$. By a general theory of Banach spaces, this embedding is isometric, see Section II.3.19 in \cite{d-sch}. Now the conclusion is evident. \qed 

\medskip

To use some results of \cite{BBK11}, we need to relate truncated Hankel and Toeplitz operators. This is the subject of the following lemma.
\begin{Lem}\label{l14}
Let $\theta$ be an inner function and let $\phi \in L^2$. Consider the truncated Hankel operator 
$\Gamma_{\phi}: \Kth \to \ov{z\Kth}$ and the truncated Toeplitz operator ${A_{\theta \phi}: \Kth \to \Kth}$. For 
all $f \in \Kth \cap L^{\infty}$ we have $\Gamma_\phi f = \bar \theta A_{\theta\phi}f$. 
\end{Lem}
\beginpf As we have seen in the proof of Proposition \ref{p2}, $P_\theta g = \theta P_{-} (\bar \theta g) - P_{-}g$ for all $g \in L^2$. Analogously, one can check that $P_{\bar \theta} g = P_{-}g - \bar \theta P_{-} (\theta g)$. Using this two formulas, we obtain
$$
\bar \theta A_{\theta \phi}f = \bar\theta\bigl(\theta P_{-}(\bar\theta\theta\phi f) - P_{-}(\theta\phi f)\bigr)
 = P_{-}(\phi f) - \bar \theta P_{-}(\theta\phi f) = \Gamma_{\phi}f
$$
for all $f \in \Kth \cap L^\infty$, as required. \qed

\medskip

\begin{Lem}\label{l11}
Let $\theta$ be an inner function. Then $\ft = (\ov{\theta H^2} + H^2) \cap L^{\infty}$. Consequently, for $\phi \in L^{\infty}$ the truncated Hankel operator $\Gamma_{\phi}: \Kth \to \ov{z\Kth}$ is the zero operator if and only if $\phi$ belongs to $\fts = (\ov{\theta^2 H^2} + H^2) \cap L^{\infty}$.
\end{Lem}
\beginpf We first check that $\mathcal{F} = (\ov{\theta H^2} + H^2) \cap L^{\infty}$ is closed in the $w^*$-topology of the space $L^\infty$. Assume that $\phi_n \in \mathcal{F}$ are such that $\lim \int \phi_n f\,dm = \int \phi f\,dm$ for every $f \in L^1$. Then $\phi \in L^\infty$ and for every function $f \in \Kth \cap zH^2$ we have $\int \phi f \,dm = 0$ because $\Kth \subset L^1$ and $\int \phi_n f\,dm = 0$ for all $n$. It follows that $\phi$ belongs to the orthogonal complement of $\ov{\Kth \cap zH^2}$ in the space $L^2$, that is, $\phi \in \mathcal{F}$. The same argument shows that $f \in L^1$ is such that $\int f \phi\,dm = 0$ for all $\phi \in \ft$ if and only if $\int f \phi\,dm = 0$ for all $\phi \in \mathcal{F}$. Since $\mathcal{F}$ and $\ft$ are $w^*$-closed, we have $\mathcal{F} = \ft$. The second part of the statement is a direct consequence of Theorem 3.1 in \cite{Sar07} and Lemma \ref{l14}. \qed

\medskip
 
\begin{Lem}\label{l12}
Let $\theta$ be an inner function. Truncated Hankel operators with continuous symbols are compact and norm-dense in the set of all compact truncated Hankel operators on $\Kth$.
\end{Lem}
\beginpf By Lemma \ref{l14}, it is sufficient to show that truncated Toeplitz operators with symbols in $\theta\ct$ are compact and norm-dense in the set of all compact truncated Toeplitz operators on $\Kth$. The first claim 
(compactness) follows from Proposition~\ref{p2}. For $\lambda \in \D$ denote $\kl = \frac{1-\ov{\theta(\lambda)}\theta}{1 - \bar \lambda z}$ and $\klt = \frac{\theta - \theta(\lambda)}{z - \lambda}$. Then $\kl, \klt \in \Kth$ and the rank-one operator $T_\lambda: h \mapsto (h, \kl)\klt$ is the truncated Toeplitz operator with symbol $\phi_{\lambda} = \frac{\theta}{z - \lambda}$, see Section 5 in \cite{Sar07}. Moreover, by Corollary 5.1 in \cite{BBK11} the linear span of the set $\{T_{\lambda}, \; \lambda \in \D\}$ is norm-dense in the set of all compact truncated Toeplitz operators on $\Kth$. Since $\phi_\lambda \in \theta \ct$ for each $\lambda \in \D$, the lemma is proved. \qed

\medskip

Next lemma is a consequence of Corollary 2.5 in \cite{BBK11}.

\begin{Lem}\label{l9}
Let $\theta$ be a one-component inner function. Then every bounded truncated Hankel operator $\Gamma$ on $\Kth$ has a bounded symbol $\phi \in L^\infty$. Moreover, there exists a constant $c_\theta$ depending only on $\theta$ such that 
$\|\Gamma\| \le \dist_{L^\infty}(\phi, \fts) \le c_\theta\|\Gamma\|$
for every symbol $\phi \in L^{\infty}$ of $\Gamma$. 
\end{Lem}
\beginpf Let $\Hth$ be the set of all bounded truncated Hankel operators on $\Kth$. It follows from Theorem 4.2 in \cite{Sar07} and Lemma \ref{l14} that $\Hth$ is closed in the weak operator topology. In particular, $\Hth$ is the Banach space with respect to the operator norm. By Corollary 2.5 of \cite{BBK11} and Lemma \ref{l14}, every bounded truncated Hankel operator $\Gamma$ on $\Kth$ has a bounded symbol $\phi_\Gamma$. The set of all bounded symbols of $\Gamma$ equals $\phi_\Gamma + \fts$ by Lemma \ref{l11}. Hence the linear mapping $\Gamma \mapsto [\phi_\Gamma]$ from $\Hth$ to $L^{\infty}/\fts$ is correctly defined and bounded by the closed graph theorem. It follows that there exists a constant $c_\theta$ depending only on $\theta$ such that $\dist_{L^\infty}(\phi,  \fts) \le c_\theta\|\Gamma\|$ for every bounded symbol $\phi$ of $\Gamma$. On the other hand, we have $\|\Gamma\| \le \dist_{L^\infty}(\phi, \fts)$ because $\|\Gamma_{\phi+\psi}\| \le \|\phi+\psi\|_{L^{\infty}}$ and $\Gamma_\psi = 0$ for every $\psi \in \fts$. \qed 

\medskip

\noindent{\bf Proof of Theorem \ref{t2}.} Let $\theta$ be a one-component inner function and let $\Gamma$ be a compact truncated Hankel operator from $\Kth$ to $\ov{z\Kth}$. Fix a number $\eps> 0$ and put $c_\eps = (1+\eps)$. By Lemma \ref{l12}, one can find compact truncated Hankel operators $\Gamma_{\phi_k}$ with symbols $\phi_k \in \ct$ such that 
$$
\Gamma = \sum_k \Gamma_{\phi_k}, \quad  \sum\|\Gamma_{\phi_k}\| \le c_\eps\|\Gamma\|. 
$$
It follows from Lemma \ref{l5} for the inner function $\theta^2$ that there exist functions $\tilde{\phi}_k \in \ct$ such that $\phi_k - \tilde{\phi}_k \in \fts$ and 
$\|\tilde{\phi}_{k}\|_{L^\infty} \le  c_\eps\dist_{L^{\infty}}(\phi_k, \fts)$. By Lemma \ref{l9}, we have 
\begin{equation}\label{eq19}
\sum_k\|\tilde{\phi}_{k}\|_{L^\infty} \le  c_\eps \sum_k \dist_{L^{\infty}}(\phi_k, \theta\fts) 
\le c_\eps c_\theta \sum\|\Gamma_{\phi_k}\| \le c_\eps^2 c_\theta \|\Gamma\|.
\end{equation}
Consider the function $\tilde\phi = \sum_k\tilde{\phi}_{k}$. By the construction, $\tilde \phi$ is continuous on $\T$ and we have $\Gamma = \sum \Gamma_{\phi_k} = \sum \Gamma_{\tilde{\phi}_k} = \Gamma_{\tilde{\phi}}$, that is, $\tilde\phi$ is the symbol of the operator $\Gamma$. Clearly, we have  $\|\Gamma\| \le \|\tilde\phi\|_{L^\infty}$. On the other hand, $\|\tilde{\phi}\|_{L^\infty} \le  c_\eps^2 c_\theta \|\Gamma\|$ by 
estimate~\eqref{eq19}. The theorem is proved. \qed

\medskip

\noindent{\bf Example.} Let $\theta$ be a finite Blaschke product. Then one can obtain Lemma~\ref{l5} (and hence Theorem \ref{t2}) using an elementary convolution argument. Indeed, take a function $\phi \in \ct$ and consider $f \in \ft$ such that 
$$
\|\phi - f\|_{L^{\infty}} = \dist_{L^\infty}(\phi, \ft).
$$
Since the subspace $\ft$ has finite codimension in $L^{\infty}$, we have $f = f_1 + \ov{\theta f_2}$ for some functions $f_1, f_2 \in H^\infty$. For $g \in L^{\infty}$ and $0<r<1$ denote by $g_{r}$ the continuous function 
$z \mapsto \hat{g}(rz)$ in the closed unit disk. Note that $\|g_r\| \le \|g\|$  because $g_r$ is the convolution of $g$ and the Poisson kernel $\frac{1-r^2}{|1 - r\zeta|^{2}}$ which is of unit norm in~$L^{1}(\T)$. Hence, 
$$
\|\phi_r - f_r\|_{L^{\infty}} \le \|\phi - f\|_{L^{\infty}}.
$$
We have $\lim_{r}\|\phi - \phi_r\|_{L^\infty} = 0$ and $\lim_{r}\|\theta - \theta_r\|_{L^\infty} =0$ because  the functions $\phi$, $\theta$ are continuous on the unit circle $\T$. Put $\tilde{f}_r = f_{1r} + \ov{\theta f_{2r}}$. Then we have
$$
\dist_{L^\infty}(\phi, \ft \cap \ct) \le \lim_{r \to 1} \|\phi - \tilde{f}_r\|_{L^\infty} = \lim_{r \to 1} \|\phi_r - f_r\|_{L^\infty} \le \|\phi - f\|_{L^{\infty}}.
$$
This shows that $\dist_{L^\infty}(\phi, \ft \cap \ct) \le \dist_{L^\infty}(\phi, \ft)$. The opposite inequality is obvious. \qed

\medskip

\noindent{\bf Remark.} The main problem with the above proof in the general case is that the equality $\lim_{r \to 1}\|\theta - \theta_r\|_{L^\infty} =0$ is false if $\theta$ is not a finite Blaschke product. 

\section{Truncated Toeplitz and Hankel operators in Schatten classes}\label{s4} 
Let $H_1$, $H_2$ be separable Hilbert spaces, and let $0 < p < \infty$. Recall that a compact operator $T: H_1 \to H_2$ belongs to the Schatten class $S^p = S^p(H_1, H_2)$ if the quantity
$$\|T\|_{S^p} = \Bigl(\sum s_k(T)^p\Bigr)^{1/p}$$
is finite, where $s_k(T)$ are the singular values of $T$ (that is, $s_k(T)$ are the eigenvalues of the compact selfajoint operator $|T|$). 
\medskip

In Section \ref{s41} we introduce some singular integral operators closely related to truncated Hankel operators. In Section \ref{s42} we collect several characterizations of compact truncated Hankel operators, extending Theorem \ref{t2}. Next in Section~\ref{s43} we present a conjecture related to truncated Toeplitz and Hankel operators in Schatten classes $S^{p}$, $0 < p < \infty$. 

\medskip

\subsection{A class of singular integral operators}\label{s41} 
Fix a complex number $\alpha$ such that $|\alpha| = 1$. Let $\sa$ denote the corresponding Clark measure of $\theta$, that is, the positive measure on the unit circle $\T$ such that
\begin{equation}\label{eq17}
\Re\left(\frac{\alpha+\theta(z)}{\alpha - \theta(z)}\right) = \int_{\T} \frac{1 - |z|^2}{|1 - \bar \xi z|^2}\,d\sa(\xi), \quad z \in \D.
\end{equation}
For a review of the theory of Clark measures see \cite{PoltSar06}. Directly from \eqref{eq17} one can see that the subset of~$\T$ where the angular boundary values of $\theta$ equals $\alpha$ is of full measure~$\sa$. In particular, the measures $\{\sa\}_{|\alpha| = 1}$ are mutually singular and singular with respect to the Lebesgue measure $m$ on $\T$. A.~G.~Pol\-torat\-ski \cite{Polt93} proved that every function $F \in \Kth$ has angular boundary values $f \in L^2(\sa)$ $\sa$-almost everywhere on the unit circle $\T$.  The function $F$ can be recovered from its trace by  
\begin{equation}\label{eq18}
F(z) = \int_{\T} f(\xi) \frac{1 - \ov{\alpha} \theta(z)}{1 - \bar\xi z}\,d\sa(\xi), \quad z \in \D.
\end{equation}
It was proved in pioneering work \cite{Cl72} by D.~N.~Clark that the mapping $V_\alpha: F \mapsto f$ is the unitary operator from $\Kth$ onto the space $L^2(\sa)$. Formula~\eqref{eq18} allows us to identify functions in $\Kth$ and their traces in $L^2(\sa)$. 

\medskip

Consider the unitary operator $H_\alpha = V_{-\alpha}V_{\alpha}^{-1}$ acting from $L^2(\sa)$ to $L^2(\sab)$. From formula \eqref{eq18} we see that for every $f\in L^2(\sa)$ and $g \in L^2(\sab)$ with separated supports we have 
\begin{equation}\label{eq21}
(H_\alpha f, g)_{L^2(\sab)} = 2\int_{\T} \left(\int_{\T}\frac{f(\xi)}{1 - \bar \xi \zeta}\,d\sa(\xi)\right) \ov{g(\zeta)}\,d\sab(\zeta),  
\end{equation}
Thus, $H_\alpha$ is the unitary Hilbert transform from $L^2(\sa)$ to $L^2(\sab)$; in the discrete setting it was characterized by Yu. Belov, T. Mengestie, and K. Seip  in \cite{Bel10}.  

\medskip

We will need the following technical lemma.   
\medskip

\begin{Lem}\label{l16}
Let $\Omega: L^2(\sa) \times L^2(\sab) \to \C$ be a sesquilinear form defined on
pairs of Lipschitz functions $f \in L^2(\sa)$, $g \in L^2(\sab)$ with $\dist(\supp f, \supp g) >0$. Assume that 
$|\Omega(f, g)| \le c\|f\|_{L^2(\sa)}\cdot\|g\|_{L^2(\sab)}$ for all such $f$, $g$ and some $c\ge0$.
Then there exists the unique linear bounded operator $T: L^2(\sa) \to  L^2(\sab)$ such that $\Omega(f, g) = (T f, g)_{L^2(\sab)}$ for all pairs $f, g$ from the domain of definition of $\Omega$. 
\end{Lem}
\beginpf Since the measures $\sa$, $\sab$ are mutually singular, for every pair of functions $f \in L^2(\sa)$, $g \in L^2(\sab)$ one can find Lipschitz functions $f_n \in L^2(\sa)$, $g_n \in L^2(\sab)$ with separated supports and such that $\lim\|f - f_n\|_{L^2(\sa)} = 0$, $\lim\|g - g_n\|_{L^2(\sab)} = 0$. This observation shows that the
domain of definition $\dmn\Omega$ of $\Omega$ is dense in $L^2(\sa) \times L^2(\sab)$. Also it is easy to check that $\Omega$ is continuous on $\dmn\Omega$. Now we can extend $\Omega$ to the whole space $L^2(\sa) \times L^2(\sab)$ by continuity and find a linear bounded operator $T: L^2(\sa) \to  L^2(\sab)$ such that $\Omega(f, g) = (T f, g)_{L^2(\sab)}$. \qed

\medskip

Lemma \ref{l16} shows that the operator $H_\alpha$ is completely determined by its sesquilinear form \eqref{eq21}. The same is true for the singular operators $C_\phi$ that will be defined below.

\medskip

Consider the measure $\nua = (\sa + \sab)/2$. Observe that $\nua$ is the Clark measure of the inner function $\theta^2$: 
$$
\Re\left(\frac{\alpha^2+\theta^2(z)}{\alpha^2 - \theta^2(z)}\right) = \int_{\T} \frac{1 - |z|^2}{|1 - \bar \xi z|^2}\,d\nua(\xi), \quad z \in \D.
$$ 
Take a function $\phi \in L^2(\nua)$. Let $M_{\phi}$ be the densely defined operator of multiplication by $\phi$ on $L^2(\nua)$. Consider the commutator 
$C_{\phi} = H_\alpha M_{\phi} - M_{\phi} H_\alpha$ as an operator acting from $L^2(\sa)$ to $L^2(\sab)$. More precisely, for every pair of functions  with separated supports, $f \in L^\infty(\sa)$, $g \in L^\infty(\sab)$, define
$$
(C_{\phi} f, g)_{L^2(\sab)} = \int_{\T}\left(\int_{\T} \frac{\phi(\xi) - \phi(z)}{1 - \bar \xi z}f(\xi)
\,d\sa(\xi)\right) \ov{g(\zeta)} \,d\sab(\zeta).
$$
Finally, introduce the unitary operator $\tilde{V}_{-\alpha}: \ov{z\Kth} \to L^2(\sab)$ which takes a function $\ov{zf} \in \ov{z\Kth}$ into $\ov{z\cdot V_{-\alpha} f} \in L^2(\sab)$. Next lemma shows that the singular operators $C_\phi$ are unitarily equivalent to the truncated Hankel operators on $\Kth$. 
\begin{Lem}\label{l15}
Let $\theta$ be an inner function and let $\ov{\phi} \in \Kthtwo$. Consider the operators $\Gamma_\phi: \Kth \to \ov{z \Kth}$ and $C_\phi: L^2(\sa) \to L^2(\sab)$. We have $\Gamma_\phi = \tilde V_{-\alpha}^{-1} C_\phi V_\alpha$. In other words, for all $f$, $g$ in $\Kth$ such that $f = f_1$ in $L^2(\sa)$ and $g = g_1$ in $L^2(\sab)$ for some Lipschitz functions $f_1, g_1$ on $\T$ with separated supports we have
\begin{equation}\label{eq22}
(\Gamma_{\phi} f, \ov{zg})_{L^2} = (C_{\phi} f, \ov{zg})_{L^2(\sab)}.
\end{equation}
In particular, the operators $\Gamma_\phi$, $C_{\phi}$ are bounded (compact, of Schatten class $S^p$) or not simultaneously and their norms coincide.	
\end{Lem}
\beginpf At first, let us check that the sesquilinear forms in \eqref{eq22} are correctly defined. For this we need to show that $f \in \Kth \cap L^\infty$ and $\phi \in L^2(\nua)$. The inclusion $f \in \Kth \cap L^\infty$ follows from \eqref{eq18} and the assumption that $f = f_1$ in $L^2(\sa)$ for a Lipschitz function $f_1$ on $\T$. We have $\phi \in L^2(\nua)$ because $\nua$ is the Clark measure for the inner function $\theta^2$ and $\ov{\phi} \in \Kthtwo$. Next, since both $\ov{\phi}$ and $zfg$ belong to $\Kthtwo$, we have
$$
(\Gamma_{\phi} f, \ov{zg})_{L^2} = (\phi f, \ov{zg})_{L^2} = (zfg, \bar{\phi})_{L^2} = (zfg, \bar{\phi})_{L^2(\nua)}. 
$$
Using \eqref{eq18} and the fact that $\theta(z) = \alpha$ for $\sa$-almost all $z \in \T$ in the sense of angular boundary values, we obtain
\begin{equation}\notag
\begin{aligned}
(\xi fg, \bar{\phi})_{L^2(\sa)} &= 2\int_{\T} \xi f_1(\xi)\phi(\xi)\int_{\T}\frac{g_1(\zeta)}{1 - \bar\zeta \xi}\,d\sab(\zeta)\,d\sa(\xi)\\ 
&= 2\iint \frac{\phi(\xi)}{\bar \xi - \bar \zeta} f_1(\xi)g_1(\zeta)\,d\sab(\zeta)\,d\sa(\xi).
\end{aligned}
\end{equation}
Note that all integrals in the above formula converge absolutely because the supports of $f_1$ and $g_1$ are separated. Analogously, we have
$$(\zeta fg, \bar{\phi})_{L^2(\sab)} 
= 2\iint \frac{\phi(\zeta)}{\bar \zeta - \bar \xi} f_1(\xi)g_1(\zeta)\,d\sa(\xi)\,d\sab(\zeta).
$$
%\begin{equation}\notag
%\begin{aligned}
%(\zeta fg, \bar{\phi})_{L^2(\sab)} &= 2\int_{\T} \zeta g_1(\zeta)\phi(\zeta)\int_{\T}\frac{f_1(\xi)}{1 - \bar\xi\zeta}\,d\sa(\xi)\,d\sab(\zeta)\\ 
%&= 2\iint \frac{\phi(\zeta)}{\bar \zeta - \bar \xi} f_1(\xi)g_1(\zeta)\,d\sa(\xi)\,d\sab(\zeta).
%\end{aligned}
%\end{equation}
Summing up this two formulas, we get
$$
(zfg, \bar{\phi})_{L^2(\nua)} = 
\iint \frac{\phi(\zeta) - \phi(\xi)}{\bar \zeta - \bar \xi} f_1(\xi)g_1(\zeta)\,d\sa(\xi)\,d\sab(\zeta)
= (C_{\phi} f, \ov{\zeta g})_{L^2(\sab)},
$$
and formula \eqref{eq22} follows. The second part of the statement is a consequence of Lemma~\ref{l16}.
 \qed

\medskip

\subsection{Compact truncated Hankel operators}\label{s42}Consider a truncated Hankel operator $\Gamma_\phi: \Kth \to \ov{z\Kth}$ with symbol $\phi \in L^2$. Denote by $\phi_s$ the orthogonal projection of $\phi$ to the subspace $\ov{\Kthtwo \cap zH^2}$ of $L^2$. It is easy to check that 
$\Gamma_\phi = \Gamma_{\phi_s}$.  The function $\phi_s$ is called the standard symbol of $\Gamma_\phi$. Standard symbols of truncated Hankel operators play the same role as the anti-analytic symbols of usual Hankel operators on $H^2$. In particular, it is possible to describe compact truncated Hankel operators in terms of the mean oscillation properties of its standard symbols.

\medskip

Let $\nu$ be a measure on the unit circle $\T$. For a function $\phi \in L^1(\nu)$ and $\eps >0$ denote 
$$
M_{\eps}(\phi) = \sup\left\{\frac{1}{\nu(\Delta)} \int_{\Delta} |\phi - \langle \phi\rangle_{\Delta, \nu}|\,d\nu, \; \mbox{$\Delta$ is an arc of $\T$ with } 0< \nu(\Delta) \le \eps \right\},
$$
where $\langle \phi\rangle_{\Delta, \nu} = \frac{1}{\nu(\Delta)}\int_{\Delta}\phi\,d\nu$. Define the space $\vmo(\nu)$ of functions of vanishing mean oscillation with respect to $\nu$ by  
$\vmo(\nu) = \{\phi \in L^1(\nu): \; \lim_{\eps \to 0}M_{\eps} \phi = 0\}.$ 
We are in position to state the description of compact truncated Hankel and Toeplitz operators.
\begin{Prop}\label{p3}
Let $\theta$ be a one-component inner function and let $\Gamma_{\phi}: \Kth \to \ov{z\Kth}$ be a truncated Hankel operator with the standard symbol $\phi \in \ov{\Kthtwo \cap zH^2}$. The following assertions are equivalent:
\begin{itemize}
\item[(1)] $\Gamma_{\phi}: \Kth \to \ov{z\Kth}$ is compact;
\item[(2)] $A_{\theta\phi}: \Kth \to \Kth$ is compact;
\item[(3)] $C_{\phi}: L^2(\sa) \to L^2(\sab)$ is compact;
\item[(4)] $\lim_{n}\|\Gamma_{\phi} - \Gamma_{\phi_n}\| = 0$ for some operators $\Gamma_{\phi_n}: \Kth \to \ov{z\Kth}$ of finite rank;
\item[(5)] $\phi \in \ct + \ov{\theta^2 H^2} + H^2$;
\item[(6)] $\phi \in \vmo(\nua)$.
\end{itemize} 
\end{Prop}
\beginpf Assertions $(1)$, $(2)$, $(3)$ and $(4)$ are equivalent for all inner functions $\theta$. 
Indeed, the equivalence $(1) \Leftrightarrow (2)$ and $(1) \Leftrightarrow (3)$ was proved in Lemma \ref{l14} and Lemma \ref{l15}, correspondingly. Evidently, $(4)$ implies $(1)$. The proof of Lemma \ref{l12} shows that
$(1)$ implies $(4)$. By Theorem 3.1 in \cite{Sar07} and Lemma \ref{l14}, we have $\Gamma_{\psi} = 0$ for $\psi \in L^2$ if and only if $\psi \in \ov{\theta^2 H^2} + H^2$. From Lemma \ref{l12} we see that $(5)$ yields~$(1)$. Now assume that $\theta$ is a one-component inner function. Then $(1)$ and $(5)$ are equivalent by Theorem \ref{t2}. Equivalence $(1) \Leftrightarrow (6)$ was proved in Proposition~4.1 in \cite{BR15}. \qed

\medskip

\noindent{\bf Remark.} In paper \cite{Roch87} R. Rochberg studied some discrete singular integral operators in connection with Toeplitz and Hankel operators on the Paley-Wiener space. That operators 
$\tilde C_\phi: L^2(\sigma) \to L^2(\sigma)$ are defined by the formula 
\begin{equation}\label{eq24}
\tilde C_\phi: f \mapsto \int_{\R\setminus\{x\}} \frac{\phi(x) - \phi(y)}{x - y}f(x)
\,d\sigma(x),
\end{equation}
where $\sigma$ is the counting measure on the set of integers~$\Z$ (note that $\sigma$ is the Clark measure $\sigma_1$ of the inner function $e^{2\pi i z}$ in the upper half-plane of the complex plane). It follows from the results by V. V. Kapustin that for a general inner function $\theta$ with the discrete Clark measure $\sigma = \sa$ every bounded operator $\tilde C_\phi$ on $L^2(\sigma)$ is unitarily equivalent to the difference of a truncated Toeplitz operator $A_\psi$ on $\Kth$ and a certain wave operator related to~$A_\psi$. For more details, see Section 3 in \cite{Kap12} and Section 2 in \cite{Kap10}.

\subsection{Truncated Hankel operators in Schatten classes $S^p$}\label{s43}
Let $\nu$ be a finite measure on the unit circle $\T$ and let $f \in L^1(\nu)$. Take a positive integer $r$. For every arc~$\Delta$ of $\T$ define the mean oscillation of $f$ of order $r$ by
$$
\osc(f, \nu, \Delta, r) = \frac{1}{\nu(\Delta)}\int_\T|f - f_{\Delta, r}|\,d\nu(\xi),
$$
where $f_{\Delta, r}$ is a polynomial of degree at most $r$ such that
$$
\int f_{\Delta, r}\bar \xi^{k}\,d\nu(\xi) = 0, \qquad k = 0,1, \ldots , r.
$$
If $\nu(\Delta) = 0$, we put $\osc(f, \nu, \Delta, r) =0$. For $p \in (0, \infty)$ let $r_p$ be the integer part of the number $1/p$. Denote by $J$ the family of all dyadic subarcs of the unit circle $\T$. It is known that the classical Besov space $B_p = B_{p,p}^{0,\,1\!/\!p}(\T)$ can be defined in terms of the mean oscillation as the set of all functions $f \in L^1$ such that
$$
\|f\|_{B_p} = \left(\sum_{\Delta \in J} \osc(f, m, \Delta, r_p)^p\right)^{1/p} < \infty,
$$
where $m$ is the Lebesgue measure on $\T$. See Theorem 1 in \cite{Dor85} or Lemma 9.9 in \cite{Ricci83} for the equivalence of the above definition of $B_p$ and the classical one. By Peller's theorem, the Hankel operator~$H_\phi: H^2 \to \ov{zH^2}$ is in the class $S^p$, $0< p< \infty$, if and only if its anti-analytic symbol $P_{-}\phi$ belongs to $B_p$, see Chapter 6 in \cite{PeBook}. 

\medskip

Now let $\nu$ be a discrete measure with isolated atoms on the unit circle $\T$, and let $\suppe\nu$ be the set of all points $\xi \in \supp\nu$ such that $\nu\{\xi\} = 0$. Assume that $\nu(\suppe\nu) = 0$. For a subset $E$ of $\T$ we will denote by $\inter E$ the interior of~$E$. Let $J_\nu$ be the family of all dyadic subarcs of $\T \setminus\suppe\nu$. More precisely, a closed arc $\Delta \subset \T$ is in $J_\nu$ if the interior of $\Delta$ is contained in a connected component $I$ of the open set $\T\setminus\suppe\nu$ and there exists an integer $k \ge 0$ such that    
$$
I = \inter \bigcup_{i=1}^{2^k} \Delta_i
$$
for some rotated copies $\Delta_i$ of $\Delta$ with disjoint interiors; $\Delta_1 = \Delta$. Define the Besov space $B_p(\nu)$ as the set of all functions $f \in L^1(\nu)$ such that 
$$
\|f\|_{B_p(\nu)} = \left(\sum_{\Delta \in J_\nu} \osc(f, \nu, \Delta, r_p)^p\right)^{1/p} < \infty.
$$
Our conjecture on truncated Toeplitz and Hankel operators in Schatten classes $S^p$ reads as follows.
\begin{Conj}
Let $\theta$ be a one-component inner function and let $\Gamma_{\phi}: \Kth \to \ov{z\Kth}$ be a truncated Hankel operator with the standard symbol $\phi \in \ov{\Kthtwo \cap zH^2}$. The following assertions are equivalent for each $p \in (0, \infty)$\textup{:}
\begin{itemize}
\item[(1)] $\Gamma_{\phi}: \Kth \to \ov{z\Kth}$ is in $S^p$;
\item[(2)] $A_{\theta\phi}: \Kth \to \Kth$ is in $S^p$;
\item[(3)] $C_{\phi}: L^2(\sa) \to L^2(\sab)$ is in $S^p$;
\item[(4)] $\Gamma_{\phi} = \sum_{n} a_n \Gamma_{\phi_n}$ for some rank-one operators $\Gamma_{\phi_n}: \Kth \to \ov{z\Kth}$ of unit norm and a sequence $\{a_n\} \subset \C$ such that $\sum_n |a_n|^p < \infty$;
\item[(5)] $\phi \in B_p + \ov{\theta^2 H^2} + H^2$;
\item[(6)] $\phi \in B_p(\nua)$.
\end{itemize} 
\end{Conj}
\noindent As before, assertions $(1)$, $(2)$, and $(3)$ are equivalent for all inner functions $\theta$. It is also easy to see that each of $(4)$ and $(5)$ implies $(1)$. Consider the case where $\theta$ is the inner function of the form $\theta: z \mapsto e^{i\pi z}$ in the upper half-plane of the complex plane. In this situation $(1)$ is equivalent to 
\begin{itemize}
\item[(5$'$)] $\phi \in B_{p}[0, 2\pi]$,
\end{itemize} 
where $B_{p}[0, 2\pi]$ is a Besov class associated with the interval $[0,2\pi]$. This result is due to 
R. Rochberg for $1 \le p < \infty$ (Theorems 5.1, 5.2 in \cite{Roch87}) and to V. Peller \cite{Pe88} for $0 < p < 1$. Also, R. Rochberg \cite{Roch87} proved that for the operators 
$\tilde C_{\phi}$ in \eqref{eq24} the assertion 
\begin{itemize}
\item[(3$'$)] $\tilde C_{\phi}: L^2(\sigma) \to L^2(\sigma)$ is in $S^p$ ($1< p <\infty$);
\end{itemize} 
is equivalent to assertion $(6)$ with $\nua = \nu_1 = \sigma$. His proof does not involve Hankel operators.

\bibliographystyle{plain} 
\bibliography{bibfile}
\enddocument